\input amstex
\documentstyle{amsppt}
\magnification=\magstep1
 \hsize 13cm \vsize 18.35cm \pageno=1
\loadbold \loadmsam
    \loadmsbm
    \UseAMSsymbols
\topmatter
\NoRunningHeads
\title A note on the generalized Euler numbers and polynomials
\endtitle
\author
  Taekyun Kim
\endauthor

\abstract In this paper, we derive some interesting symmetric
properties for the generalized Euler numbers and polynomials
attached to $\chi$ using the $p$-adic invariant integral on $\Bbb
Z_p$. \endabstract
\thanks  2000 AMS Subject Classification: 11B68, 11S80
\newline  The present Research has been conducted by the research
Grant of Kwangwoon University in 2008
\endthanks
\endtopmatter

\document

{\bf\centerline {\S 1. Introduction}}

 \vskip 20pt

Let $p$ be an odd prime number. Throughout this paper $\Bbb Z_p ,$
$\Bbb Q_p ,$ $\Bbb C,$ and $\Bbb C_p$ will, respectively, denote the
ring of $p$-adic rational integers, the field of $p$-adic rational
numbers, the complex number field, and the completion of algebraic
closure of $\Bbb Q_p .$  The normalized valuation in $\Bbb C_p$ is
denoted by $|\cdot |_p$ with $|p|_p =\frac{1}{p} .$ We say that $f$
is a uniformly differentiable function at a point $a \in\Bbb Z_p $
and denote this property by $f\in UD(\Bbb Z_p )$, if the difference
quotients $F_f (x,y) = \dfrac{f(x) -f(y)}{x-y} $ have a limit
$l=f^\prime (a)$ as $(x,y) \to (a,a)$. For $f\in UD(\Bbb Z_p )$, let
us start with the expression
$$\eqalignno{ & \sum_{0\leq j < p^N} (-1)^j f(j)
=\sum_{0\leq j < p^N} f(j) \mu (j +p^N \Bbb Z_p ), }$$ representing
a $p$-adic analogue of Riemann type sums for $f$, see [4-17]. The
integral of $f$ on $\Bbb Z_p$ will be defined as limit ($n \to
\infty$) of those sums, when it exists. The $p$-adic invariant
integral on $\Bbb Z_p$ of the function $f\in UD(\Bbb Z_p )$ is
defined as
$$I(f)=\int_{\Bbb Z_p}f(x)dx=\lim_{N\rightarrow
\infty}\sum_{x=0}^{p^N-1}f(x)(-1)^x, \text{ see [4-17]}.\tag1$$

In this paper we investigate some symmetric properties related to
the $p$-adic invariant integral on $\Bbb Z_p$. By using these
properties, we derive some interesting identities related to the
generalized Euler numbers and  polynomials attached to $\chi$.

\vskip 20pt

{\bf\centerline {\S 2. Some symmetric identities of the generalized
Euler polynomials attached to $\chi$}} \vskip 10pt

The $n$-th Euler numbers $E_n$ are defined as

$$\frac{2}{e^t +1}= e^{Et}=\sum_{n=0}^{\infty} E_n \frac{t^n}{n!},
\text{ (cf.[1-13]),} \tag2$$ with the usual convention about
replacing $E^n$ by $E_n.$

The $n$-th Euler polynomials $E_n(x)$  are also defined as
$$E_n(x)=\sum_{l=0}^n \binom{n}{l}x^{n-l}E_l, \text{ where $\binom{n}{l}=\frac{n\cdot (n-1)\cdots(n-l+1)}{l!}$}, \text{ (see [1-6])}
.$$ From (1) we can
easily derive
$$I(f_1)+I(f)=2f(0), \text{ where
$f_1(x)=f(x+1)$}.\tag3$$ By continuing this process, we see that
$$I(f_n)+(-1)^{n-1}I(f)=2\sum_{l=0}^{n-1}(-1)^{n-1-l}f(l),
\text{ where $f_n(x)=f(x+n)$.}$$ When $n$ is an odd positive
integer, we obtain
$$ I(f_n)+I(f)=2\sum_{l=0}^{n-1}(-1)^lf(l).
\tag4$$ If $n\in\Bbb N$ with $n\equiv 0 $ ($\mod 2$), then we have
$$q^n I(f_n)-I(f)=2\sum_{l=0}^{n-1}(-1)^{l-1}f(l)
.\tag5$$ For a fixed odd positive integer $d$ with $(p,d)=1$, set
$$\split
& X=X_d = \lim_{\overleftarrow{n} } \Bbb Z/ dp^n \Bbb Z ,\ X_1 =
\Bbb Z_p , \cr  & X^\ast = \underset {{0<a<d p}\atop {(a,p)=1}}\to
{\cup} (a+ dp \Bbb Z_p ), \cr & a+d p^n \Bbb Z_p =\{ x\in X | x
\equiv a \pmod{dp^n}\},\endsplit$$ where $a\in \Bbb Z$ lies in
$0\leq a < d p^n$. Let $\chi$ be the Dirichlet's character with
conductor $d(=odd)\in\Bbb N$. If we take $f(x)=\chi(x) e^{tx}$, then
we have

$$\int_{X}\chi(y)e^{(x+y)t}dy
=\frac{2\sum_{l=0}^{d-1}(-1)^l\chi(l)e^{lt}}{e^{dt}+1}e^{xt}.
\tag6$$ Now we define the generalized Euler polynomials attached to
$\chi$ as follows:
$$\frac{2\sum_{l=0}^{d-1}(-1)^l\chi(l)e^{lt}}{e^{dt}+1}e^{xt}=\sum_{n=0}^{\infty}E_{n,\chi}(x)\frac{t^n}{n!}.
\tag7$$ In particular $x=0$, $E_{n,\chi}(=E_{n,\chi}(0))$ are called
the $n$-th generalized Euler numbers attached to $\chi$. From (6)
and (7), we note that
$$E_{n, \chi}(x)=\int_{X}\chi(y)(x+y)^n dy. \tag8$$
For $ n \in \Bbb N$ with $ n\equiv 1$(mod $2$),  we have
$$\int_{X}\chi(x)e^{(nd+x)t}dx+\int_{X}\chi(x)e^{xt}dx=2\sum_{l=0}^{nd-1}(-1)^l\chi(l)e^{lt}.
\tag9$$ Let
$$T_{k, \chi}(n)=2\sum_{l=0}^n(-1)^l\chi(l)l^k. \tag10$$
It is not difficult to show that
$$\int_{X}e^{(x+nd)t}\chi(x)dx
+\int_{X}\chi(x)e^{xt}dx=\frac{2\int_{X}e^{xt}\chi(x)dx}{\int_{X}e^{ndxt}dx}=\sum_{k=0}^{\infty}T_{k,\chi}(nd-1)\frac{t^k}{k!}.
\tag11$$ Let $w_1, w_2 \in \Bbb N$ with $ w_1 \equiv 1$(mod $2$),
$w_2 \equiv 1$(mod $2$). Then we consider following double integral.
$$\aligned
&\frac{\int_{X}\int_{X}e^{(w_1x_1+w_2x_2)t}\chi(x_1)\chi(x_2)dx_1
dx_2}{\int_{X}e^{dw_1w_2xt}dx}\\
&=\left(\frac{2(e^{dw_1w_2t}+1)}{(e^{w_1dt}+1)(e^{w_2dt}+1)}\right)\left(\sum_{a=0}^{d-1}\chi(a)e^{w_1at}(-1)^a\right)
\left(\sum_{b=0}^{d-1}\chi(b)e^{w_2bt}(-1)^b \right).
\endaligned\tag12$$
By (8) and (11), we see that
$$\int_X \chi(x)(dn+x)^k dx +\int_x \chi(x)x^k dx =T_{k,
\chi}(nd-1).\tag13$$ That is,
$$E_{k,x}(nd)+E_{k,\chi}=T_{k, \chi}(nd-1).$$
Let
$$T_{\chi}(w_1, w_2)=\frac{\int_X \int_X
\chi(x_1)\chi(x_2)e^{(w_1x_1+w_2x_2+w_1w_2x)t}dx_1 dx_2}{\int_x e^{d
w_1 x_3 t}dx_3}.\tag14$$
 From (12) and (14), we note that
$$\aligned
&T_{\chi}(w_1, w_2 )\\
&=\left(\frac{2(e^{dw_1 w_2t}+1)e^{w_1w_2xt}}{(e^{w_1dt}+1)(e^{w_2
dt}+1) }\right)\left(\sum_{a=0}^{d-1}\chi(a)(-1)^ae^{w_1 at}\right)
\left(\sum_{b=0}^{d-1}\chi(b)(-1)^b e^{w_2 bt} \right),
\endaligned\tag15$$
and $T_{\chi} (w_1, w_2)$ is symmetric in $w_1$ and $w_2$. By
(6),(8),(11),(14), and (15), we see that
$$T_{\chi}(w_1,
w_2)=\sum_{l=0}^{\infty}\left(\sum_{i=0}^l\binom{l}{i}E_{i,
\chi}(w_2 x)T_{l-i,
\chi}(dw_1-1)w_1^iw_2^{l-i}\right)\frac{t^l}{l!}.\tag15$$ By the
symmetric property of $T_{\chi}(w_1, w_2)$ in $w_1$ and $ w_2$, we
also see that
$$T_{\chi}(w_1, w_2)=\sum_{l=0}^{\infty}\left(\sum_{i=0}^l \binom{l}{i}w_2^i
w_1^{l-i}E_{i, \chi}(w_1 x) T_{l-i,
\chi}(dw_2-1)\right)\frac{t^l}{l!}. \tag16$$

By comparing the coefficients on the both sides of (15) and (16), we
obtain the following theorem.

\proclaim{ Theorem 1} Let $\chi$ be the Dirichlet's character with
an odd conductor $d\in\Bbb N$. For  $w_1, w_2, d\in \Bbb N$ with $d
\equiv 1 (\mod 2)$, $w_1\equiv 1(\mod 2)$, $w_2\equiv 1 (\mod 2)$,
we have
$$\sum_{i=0}^l \binom{l}{i}E_{i, \chi}(w_2 x) T_{l-i, \chi}(d
w_1-1)w_1^i w_2^{l-i} =\sum_{i=0}^l \binom{l}{i}E_{i, \chi}(w_1 x)
T_{l-i, \chi}(d w_2-1)w_2^i w_1^{l-i}. $$
\endproclaim

Remark. Setting $x=0$ in Theorem 1 we obtain
$$\sum_{i=0}^l \binom{l}{i}E_{i, \chi} T_{l-i, \chi}(d
w_1-1)w_1^i w_2^{l-i} =\sum_{i=0}^l \binom{l}{i}E_{i, \chi} T_{l-i,
\chi}(d w_2-1)w_2^i w_1^{l-i}. $$

From (14), we can derive
$$\aligned
T_{\chi}(w_1, w_2)&=\sum_{l=0}^{dw_1-1}(-1)^l\int_X
e^{(x_1+w_2x+\frac{w_2}{w_1}l)tw_1}\chi(x_1) dx_1\\
&=\sum_{n=0}^{\infty}\left(\sum_{l=0}^{dw_1-1}(-1)^lE_{n,\chi}(w_2x+\frac{w_2}{w_1}l)w_1^n\right)\frac{t^n}{n!}.
\endaligned\tag17$$

On the other hand,
$$T_{\chi}(w_1,
w_2)=\sum_{n=0}^{\infty}\left(\sum_{l=0}^{dw_2-1}(-1)^lE_{n,
\chi}(w_1x+\frac{w_1}{w_2}l)w_2^n\right)\frac{t^n}{n!}. \tag18$$

By (17) and (18), we obtain the following theorem.

\proclaim{ Theorem 2}  For  $w_1, w_2 \in \Bbb N$ with $w_1\equiv
1(\mod 2)$, $w_2\equiv 1 (\mod 2)$, we have
$$w_1^n\sum_{l=0}^{dw_1-1}(-1)^lE_{n,\chi}(w_2x+\frac{w_2}{w_1}l)
=w_2^n\sum_{l=0}^{dw_2-1}(-1)^lE_{n, \chi}(w_1x+\frac{w_1}{w_2}l)$$

\endproclaim

Remark. From the equations of the $p$-adic invariant integral on
$\Bbb Z_p$, some interesting identities of the Euler polynomials are
derived in [6]. In this paper, we have studied the symmetric
properties of the generalized Euler polynomials attached to $\chi$.
To derive the symmetric identities  of the generalized Euler
polynomials attached to $\chi$, we used the symmetric properties of
$p$-adic invariant integrals on $\Bbb Z_p$.

\vskip 10pt

\quad Taekyun Kim

\quad Division of General Education-Mathematics, Kwangwoon
University, Seoul

\quad 139-701, S. Korea
 e-mail:\text{ tkkim$\@$kw.ac.kr}
 \vskip 10pt

\enddocument